\documentclass[twocolumn,fleqn]{autart}
\usepackage{amsmath,amssymb}
\usepackage{graphicx}
\usepackage{color}
\definecolor{change}{rgb}{0, 0, 0}
\DeclareMathAlphabet{\mathbbold}{U}{bbold}{m}{n}

\newcommand{\R}{\mathbbold{R}}

\newcommand{\N}{\mathbbold{N}}

\begin{document}

\begin{frontmatter}

\title{\textcolor{change}{Output-Feedback} Control of the Semilinear Heat Equation\\ via the $L^2$ Residue Separation and Harmonic Inequality}

\thanks{The paper was not presented at any conference. }

\author[First]{Anton~Selivanov}\ead{a.selivanov@sheffield.ac.uk},
\author[Second]{Emilia~Fridman}\ead{emilia@tauex.tau.ac.il}

\address[First]{School of Electrical and Electronic Engineering, The University of Sheffield, UK}
\address[Second]{School of Electrical Engineering, Tel Aviv University, Israel}


\begin{abstract}\textcolor{change}{
A popular approach to designing finite-dimensional boundary controllers for partial differential equations (PDEs) is to decompose the PDE into independent modes and focus on the dominant ones while neglecting highly damped residual modes. However, the neglected modes can adversely affect the overall system performance, causing spillover. The $L^2$ residue separation method was recently introduced to eliminate spillover in the state-feedback control design. In this paper, we extend this method to finite-dimensional output-feedback control, where the output is contaminated by the residual modes. To deal with the output residue, we introduce a new harmonic inequality that optimally bounds it. We develop the approach for a 1D heat equation with unknown nonlinearity, where boundary temperature measurements are used to control heat flux at the opposite boundary. By exploiting the connection between $L^2$ residue separation and $H_\infty$ theory, we show that the class of admissible nonlinearities can only increase with higher controller order.}%
\end{abstract}

\end{frontmatter}
\section{Introduction}
Modal decomposition is a widely used technique for designing finite-dimensional controllers for systems governed by partial differential equations (PDEs). The approach involves representing the solution as a Fourier series and concentrating on a finite number of dominant modes, while disregarding the highly damped residual modes~\cite{Athans1970,Triggiani1980,Balas1982a,Curtain1993}. This way finite-dimensional control methods can be adapted to tackle infinite-dimensional PDE systems. The critical limitation of this approach is spillover: the neglected modes can adversely affect the overall system performance \cite{Meirovitch1983,Bontsema1988,Hagen2003}. Spillover has been addressed using residual mode filters \cite{Balas1988,Moheimani1998,Harkort2011}, spectral properties of linear operators~\cite{Lasiecka1983,Curtain1984,Orlov2017}, small-gain techniques \cite{Karafyllis2019b,Xia2023}, and Lyapunov functionals \cite{Karafyllis2019k,Karafyllis2019i,Karafyllis2021}. These qualitative results provided valuable insights, such as establishing stability when a sufficiently large number of modes are considered, though they do not specify the exact number of modes required. 

Achieving precise quantitative results requires a more refined analysis of the residual modes, which can be conducted using Lyapunov functionals~\cite{Hagen2003,Hagen2006,Selivanov2018e,Katz2020a,Katz2021c,Katz2022b,Lhachemi2022c}. A critical step in
such analyses is handling the interaction between the control input and the residual modes, which is typically done using Young’s inequality to separate cross terms. 

An improved approach to residual modes separation for spillover avoidance was introduced \textcolor{change}{for the state-feedback case} in~\cite{Selivanov2025a}. The key insight was that, when the controller is designed without explicitly accounting for the residual dynamics, the residual modes can be treated as being disturbed by an unknown control input. By determining the corresponding $L^2$ gains for these disturbances and summing them to infinity, a combined $L^2$ gain is obtained, which characterizes the influence of the control input on the residual modes. This $L^2$ gain can then be utilized to avoid spillover when designing a controller for the dominant modes. \textcolor{change}{The $L^2$ residue separation was refined in \cite{Selivanov2024a} for the design of guaranteed-cost state-feedback controllers for the semilinear heat equation}, demonstrating up to a 90\% cost reduction compared to methods based on Young’s inequality.

\textcolor{change}{The present paper extends the $L^2$ residue separation method to finite-dimensional \textit{output-feedback} control. The main challenge in this setting is that the measured output is contaminated by the infinite-dimensional residue, making the separation process substantially more intricate than in the state-feedback case. To overcome this, we introduce a new harmonic inequality (Lemma~\ref{lem:harmonic}) that provides an optimal bound on the residue in terms of the weighted $l^2$ norm of its Fourier coefficients. This result leads to significantly tighter residue bounds, enabling the design of low-order, spillover-free controllers.} 

\textcolor{change}{Furthermore, we exploit the connection between $L^2$ residue separation and $H_\infty$ theory to handle unknown nonlinearities and design the controller and observer gains. The link with $H_\infty$ theory enables us, for the first time, to show that the class of admissible nonlinearities does not shrink as the controller order increases (Proposition~\ref{prop:N}), and to provide a theoretical justification for zeroing certain controller and observer gains in the linear case (Remark~\ref{cor:linear}). Moreover, the proposed $L^2$ residue separation approach avoids the need for the lifting transformation \cite{Katz2022b,Kang2019}, thereby simplifying the analysis of control inputs with time-varying delays and enabling a sample-and-hold implementation (Section~\ref{sec:sampled-data}) instead of the less practical generalized hold used in previous works.} 
 

\textit{Notation:} $\N_0:=\N\cup\{0\}$, $|\cdot|$ is the Euclidean norm, $\|\cdot\|$ and $\langle \cdot,\cdot\rangle$ are the $L^2$ norm and product. The Sobolev spaces $H^1$ ($H^2$) consist of square-integrable functions whose first (and second) weak derivatives are also square-integrable. If $P$ is a symmetric matrix, $P>0$ means that it is positive definite with the symmetric elements sometimes denoted by ``$*$''. If $P>0$, then $|z|_P^2:=z^\top Pz$. The spectral radius of $A\in\R^{n\times n}$ is $\rho(A)=\max\{|\lambda_1|,\ldots,|\lambda_n|\}$, where $\{\lambda_i\}$ are the eigenvalues of $A$. Partial derivatives are denoted by indices, e.g., $z_x=\frac{\partial z}{\partial x}$. 
\subsection{The harmonic inequality}
\textcolor{change}{The following inequality is a key element facilitating the $L^2$ residue separation in the output-feedback setting.}

\begin{lem}[Harmonic inequality]\label{lem:harmonic}~\\
\begin{equation}\label{harmonicGeneral}
\left(\sum_{n=1}^\infty z_n\right)^2\le\sum_{n=1}^\infty\mu_nz_n^2
\end{equation}
for all $(z_n)_{\N}\in l^1$ if and only if $\mu_n$ satisfy 
\begin{equation*}
\sum_{n=1}^\infty\mu_n^{-1}\le1,\quad \mu_n>0 \quad\text{(harmonic condition).}
\end{equation*}
\end{lem}
{\bf Proof.} For each $M\in\N$, define $\mathbf{z}=\left[z_{1},\ldots,z_{M}\right]^\top$, $\mathbf{1}=\left[1,\ldots,1\right]^\top\in\R^{M}$, and $D=\operatorname{diag}\{\mu_{1},\dots,\mu_{M}\}$. Then \eqref{harmonicGeneral} holds if and only if 
\begin{equation*}
\textstyle\mathbf{z}^\top\mathbf{1}\mathbf{1}^\top\mathbf{z}=\left(\sum_{n=1}^{M}z_n\right)^2\le\sum_{n=1}^{M}\mu_nz_n^2=\mathbf{z}^\top D\mathbf{z},\:\textcolor{change}{\forall M\in\N.}
\end{equation*}
By the Schur complement lemma, $\mathbf{1}\mathbf{1}^\top\le D$ is equivalent~to 
\begin{equation*}
\textstyle\left[\begin{smallmatrix}
    D & \mathbf{1} \\ \mathbf{1}^\top & 1
\end{smallmatrix}\right]\ge0\quad\iff\quad 
\mathbf{1}^\top D^{-1}\mathbf{1}=\sum_{n=1}^{M}\mu_n^{-1}\le 1, 
\end{equation*}
which is equivalent to the harmonic condition \textcolor{change}{as $M\to\infty$}. 

\begin{rem}[Jensen's inequality]\label{lem:Jensen}
The sufficiency of the harmonic condition also follows from Jensen's inequality: 
\begin{equation*}
\begin{array}{rcl}
\left(\sum_{n=1}^{\infty}z_n\right)^2&=&\left(\sum_{n=1}^{\infty}\mu_n^{-1}\left(\mu_nz_n\right)\right)^2\\
&\le&\sum_{n=1}^{\infty}\mu_n^{-1}\left(\mu_nz_n\right)^2=\sum_{n=1}^{\infty}\mu_nz_n^2. 
\end{array}
\end{equation*}
\end{rem}
\section{Modal decomposition of the semilinear heat PDE with boundary input and output}\label{sec:modalDecomposition}
Consider the semilinear heat equation: 
\begin{subequations}\label{HeatBVP}
\begin{align}
& z_t(x,t)=z_{xx}(x,t)+qz(x,t)+f(x,t,z(\cdot,t)),\label{PDE}\\
& z_x(0,t)=0, \quad z_x(\pi, t)=u(t),\label{BC}\\
& y(t)=z(0,t)\label{Output}
\end{align}
\end{subequations}
with the state $z\colon[0,\pi]\times[0,\infty)\to\R$, control input $u\colon[0,\infty)\to\R$, measured output $y\colon[0,\infty)\to\R$, reaction coefficient $q>0$, and \textcolor{change}{nonlinear $f$ satisfying }
\begin{subequations}\label{fSector}
\textcolor{change}{\begin{align}
& f\in C([0,\infty)\times L^2(0,\pi);L^2(0,\pi)),\\
& f(t,0)=0,\quad\forall t\ge0,\\
& \exists\sigma>0\colon \forall t\ge0, z_1,z_2\in L^2(0,\pi),\nonumber\\
&\hspace{1cm}\|f(t,z_1)-f(t,z_2)\|\le\sigma\|z_1-z_2\|. 
\end{align}}%
\end{subequations}
We assume that $f$ is unknown, but $\sigma$ is known. If the diffusion coefficient in front of $z_{xx}$ is not $1$, or the spatial domain is not $[0,\pi]$, the equation can be transformed into the form \eqref{HeatBVP} using the change of variables $\tilde z(x,t)=z(ax-x_0,bt)$ with suitable $a$, $b$, and~$x_0$. 
Note that the reaction term, $qz$, can be incorporated into $f$, which would increase the Lipschitz constant $\sigma$, but we treat it separately to obtain more precise conditions. 

Our goal is to design a finite-dimensional output-feedback controller that stabilizes \eqref{HeatBVP} under the Lipschitz condition \eqref{fSector}. To achieve this, we perform modal decomposition, which represents the solution of \eqref{HeatBVP} as a sum of modes that are independent in the linear case. This is achieved by diagonalizing the operator 
\begin{equation*}
\mathcal{A}\varphi=-\varphi'',\: D(\mathcal{A})=\{\varphi\in H^2(0,\pi)\mid \varphi'(0)=0=\varphi'(\pi)\}. 
\end{equation*}
\textcolor{change}{Its eigenvalues and eigenfunctions 
\begin{equation*}
\begin{aligned}
&\lambda_n=n^2,\qquad n\in\N_0,\\
&\varphi_n(x)=\left\{\begin{aligned}
&1/\sqrt{\pi},&&n=0,\\
&\sqrt{2/\pi}\cos nx,&&n\in\N, 
\end{aligned}\right.
\end{aligned}
\end{equation*}
form an orthonormal basis of $L^2(0,\pi)$.} Therefore, the state can be presented as the Fourier series, 
\begin{equation*}
\textstyle z(\cdot,t)\stackrel{L^2}{=}\sum_{n=0}^{\infty} z_n(t) \varphi_n(\cdot),\qquad z_n(t):=\langle z(\cdot,t), \varphi_n\rangle. 
\end{equation*}
The Fourier coefficients, $z_n(t)$, satisfy 
\begin{equation*}
\begin{aligned}
\dot z_n(t)&\!=\!\langle z_t(\cdot,t), \varphi_n\rangle\\
&\!\!\!\stackrel{\eqref{PDE}}{=}\!
\langle z_{xx}(\cdot,t),\varphi_n\rangle+q\langle z(\cdot,t),\varphi_n\rangle+\langle f(\cdot,t,z(\cdot,t)),\varphi_n\rangle.
\end{aligned}
\end{equation*}
Since $\varphi_n\in D(\mathcal{A})$ and $\varphi_n''=-\lambda_n\varphi_n$, integrating by parts twice, we obtain 
\begin{equation*}
\begin{aligned}
\langle z_{xx}(\cdot,t),\varphi_n\rangle&\stackrel{\phantom{\eqref{BC}}}{=}\left[z_x(\cdot,t)\varphi_n\right]_0^\pi-\left[z(\cdot,t)\varphi_n'\right]_0^\pi+\langle z(\cdot,t),\varphi_n''\rangle\\
&\stackrel{\eqref{BC}}{=}\varphi_n(\pi)u(t)-\lambda_nz_n(t). 
\end{aligned}
\end{equation*}
We will show that $z(\cdot,t)\in H^1(0,\pi)$, which implies 
\begin{equation*}
\textstyle y(t)=z(0,t)=\sum_{n=0}^\infty z_n(t)\varphi_n(0). 
\end{equation*}
Therefore, the modal decomposition of \eqref{HeatBVP} is 
\begin{equation}\label{modalDecomposition} 
\begin{aligned}
\dot{z}_n(t)&=-(\lambda_n-q)z_n(t)+b_nu(t)+f_n(t),\quad n\in\N_0, \\
y(t)&\textstyle=\sum_{n=0}^\infty c_nz_n(t), 
\end{aligned}
\end{equation}
where 
\begin{equation*}
\begin{aligned}
&f_n(t)=\langle f\left(\cdot,t,z(\cdot,t)\right),\varphi_n\rangle,\qquad n\in\N_0,\\
&b_n=\varphi_n(\pi)=\left\{
\begin{aligned}
&1/\sqrt{\pi},&&n=0,\\
&(-1)^n\sqrt{2/\pi},&&n\in\N,
\end{aligned}\right.\\
&c_n=\varphi_n(0)=\left\{
    \begin{aligned}
    &1/\sqrt{\pi},&&n=0,\\
    &\sqrt{2/\pi},&&n\in\N.
    \end{aligned}
    \right.
\end{aligned}
\end{equation*}
The modal decomposition \eqref{modalDecomposition} comprises infinitely many nonlinear systems with the linear part $-(\lambda_n-q)\to-\infty$ as $n\to\infty$. We will design a finite-dimensional output-feedback controller using the first $N$ modes ($n=0,1,\ldots,N-1$) with $N$ such that 
\begin{equation}\label{N}
    \lambda_N=N^2>q+\sigma. 
\end{equation}
That is, we consider all the unstable modes in the control design. Separating the first $N$ modes in \eqref{modalDecomposition}, we obtain 
\begin{subequations}\label{HeatModalDecomposition}
\begin{align}
\dot{z}^N(t)&=Az^N(t)+Bu(t)+F(t),\label{Nmodes}\\
\dot{z}_n(t)&=-\left(\lambda_n-q\right) z_n+b_nu(t)+f_n(t), \quad n\ge N, \label{residue}\\
y(t)&=Cz^N(t)+\zeta(t), \label{output}
\end{align}
\end{subequations}
where 
\begin{equation}\label{notations}
\begin{aligned}
&z^N=[z_0,\ldots,z_{N-1}]^\top , \quad F=[f_0,\ldots,f_{N-1}]^\top,\\
&A=\operatorname{diag}\{q-\lambda_0,\ldots,q-\lambda_{N-1}\},\\
&B=[b_0,b_1,\ldots,b_{N-1}]^\top,\\
&C=[c_0,c_1,\ldots,c_{N-1}],\\
&\textstyle\zeta(t)=\sum_{n=N}^\infty c_nz_n(t)=\sqrt{2/\pi}\sum_{n=N}^\infty z_n(t). 
\end{aligned}
\end{equation}
If the system is linear, i.e., $F\equiv0$ and all $f_n=0$, a linear \textit{state}-feedback that stabilizes \eqref{Nmodes} will also stabilize the full system \eqref{HeatModalDecomposition} and, therefore, \eqref{HeatBVP}. Such truncation works in this simplified case because stable residual modes \eqref{residue} do not affect \eqref{Nmodes} and remain stable under the exponentially vanishing external input \textcolor{change}{\cite{Triggiani1980}}. In the presence of nonlinearities or \textit{output}-feedback, such truncation may lead to spillover because the residual modes \eqref{residue} affect \eqref{Nmodes} either through the nonlinearity or the output feedback that uses~\eqref{output}. The next section proposes a method to avoid spillover by accounting for the residual dynamics \eqref{residue} when designing a controller for \eqref{Nmodes}. 
\section{Controller design and stability analysis\\ via the $L^2$ residue separation}\label{sed:main}
\subsection{Dynamic controller structure}
The structure of the finite-dimensional dynamic controller is selected based on the dynamics of the first $N$ modes \eqref{Nmodes}. The controller parameters are then selected to avoid the spillover from the residual dynamics \eqref{residue}. 

Since all the eigenvalues, $\lambda_n$, are different, the pair $(A,B)$ is controllable, e.g., by the Hautus lemma \cite[Lemma~3.3.1]{Sontag2013}. That is, $A-BK$ is stable for some $K\in\R^{1\times N}$. Similarly, $(A,C)$ is observable, meaning that $A-LC$ is stable for some $L\in\R^N$. This motivates the following observer-based controller 
\begin{subequations}\label{control}
\begin{align}
\dot{\hat z}^N(t)&=(A\!+\!\gamma\sigma X)\hat z^N(t)\!+\!Bu(t)\!-\!L(C\hat z^N(t)\!-\!y(t)),\label{observer}\\
u(t)&=-K\hat z^N(t)\label{u}
\end{align}
\end{subequations}
with $X\in\R^{N\times N}$ and $\gamma\in(0,\infty)$ defined below, and $\hat z^N(0)=0$. \textcolor{change}{Since $F$ is unknown, we treat it as an external disturbance. Namely, following the $H_\infty$ approach, we account for the worst-case $F=\gamma\sigma Xz^N$ (see (6.2.22) in \cite{Green2012}) by modifying the linear part of \eqref{observer} as in \cite[Section~8.2.2]{Green2012}}. 

The well-posedness of \eqref{HeatBVP} under \eqref{control} can be established in a manner similar to \cite[Section~2.2]{Katz2022b}: If $z(\cdot,0)\in L^2(0,\pi)$, the existence of a unique mild solution $z\in C([0,\infty),L^2)$ follows from \cite[Theorem~6.1.2]{Pazy1983}. If $z(\cdot,0)\in H^1(0,\pi)$, by \cite[Theorem~6.3.1]{Pazy1983}, it becomes the unique classical solution 
\begin{equation*}
\begin{aligned}
    &z\in C([0,\infty),L^2)\cap C^1((0,\infty),L^2),\\
    &z(\cdot,t)\in H^2(0,\pi),\quad\forall t>0. 
\end{aligned}
\end{equation*}
In particular, $z(\cdot,t)\in H^1(0,\pi)$, which we used to obtain the modal decomposition \eqref{modalDecomposition}. 
\subsection{The $L^2$ residue separation}
Introduce the estimation error of the dominating modes
\begin{equation*}
    e^N(t)=\hat z^N(t)-z^N(t). 
\end{equation*}
The innovation term in \eqref{observer} can be expressed as 
\begin{equation*}
    C\hat z^N(t)-y(t)=Ce^N(t)-\zeta(t)
\end{equation*}
with $C$ and residue $\zeta(t)$ defined in \eqref{notations}. Taking the difference between \eqref{observer} and \eqref{Nmodes}, we obtain 
\begin{equation}\label{e}
    \dot e^N=(A+\gamma\sigma X-LC)e^N+L\zeta-\tilde F, 
\end{equation}
where $\tilde F=F-\gamma\sigma Xz^N$ can be interpreted as the deviation from the worst-case nonlinearity. The stability of \eqref{HeatModalDecomposition}, \eqref{control} follows from that of \eqref{HeatModalDecomposition}, \eqref{e} \textcolor{change}{with $u=-K(z^N+e^N)$}. To study the stability of \eqref{HeatModalDecomposition}, \eqref{e}, we employ the following Lyapunov functional 
\begin{equation}\label{LyapunovFunctional}
\textstyle V=V_z+V_e+V_\infty,
\end{equation}
where 
\begin{equation*}
\begin{aligned}
&V_z=|z^N|_X^2,&& 0<X\in\R^{N\times N},\\
&V_e=|e^N|_Y^2,&& 0<Y\in\R^{N\times N},\\
&V_\infty\textstyle=\gamma^{-1}\sum_{n=N}^{\infty} z_n^2,&& 0<\gamma\in(0,\infty). 
\end{aligned}
\end{equation*}
To deal with the nonlinearity in the stability analysis, we will use the S-procedure~\cite{Yakubovic1977}, which leverages a quadratic constraint on the nonlinearity. Namely, by Parseval's theorem, \eqref{fSector} implies 
\begin{equation}\label{SProc_f}
\textstyle 0 \le\frac{\sigma}{\gamma}\sum_{n=0}^{\infty}z_n^2(t)-\frac1{\gamma\sigma}\sum_{n=0}^{\infty} f_n^2(t). 
\end{equation}
The right-hand side of \eqref{SProc_f} is added to $\dot V$ to introduce the negative terms $-f_n^2$ compensating the cross terms with~$f_n$. A similar method for addressing the nonlinearity was employed in \cite{Katz2022b}, where the nonlinearity was assumed to be known. \textcolor{change}{To treat the unknown nonlinearity, we combine the S-procedure with the $H_\infty$ framework \cite[Section 8.2.2]{Green2012} by decomposing $F$ into the worst-case component $\gamma\sigma Xz^N$ and the deviation $\tilde F$, given below \eqref{e}.}

The main idea of the $L^2$ residue separation is to decompose the derivative of the Lyapunov functional as follows 
\begin{subequations}\label{Vdot}
\begin{align}
&\dot V\le\dot V_z+\dot V_e+\dot V_\infty\notag\\
&\textstyle\hspace{.5cm}\pm u^2\pm\gamma^{-2}\zeta^2\pm|B^\top Xe^N|^2\pm\frac1{\gamma\sigma}|\tilde F|^2+\eqref{SProc_f}\notag\\
&\textstyle{}=\left[\dot V_z\!+\!u^2\!-\!|B^\top\! Xe^N|^2\!+\!\frac{\sigma}{\gamma}|z^N|^2\!-\!\frac1{\gamma\sigma}|F|^2\!+\!\frac1{\gamma\sigma}|\tilde F|^2\right]\label{Vzdot}\\
&\textstyle{}+\left[\dot V_e-\gamma^{-2}\zeta^2+|B^\top Xe^N|^2-\frac1{\gamma\sigma}|\tilde F|^2\right]\label{Vedot}\\
&\textstyle{}+\Bigl[\dot V_\infty+\!\gamma^{-2}\zeta^2\!-\!u^2\!+\frac{\sigma}{\gamma}\sum_{n=N}^\infty z_n^2-\frac{1}{\gamma\sigma}\sum_{n=N}^\infty f_n^2\Bigr]\label{Vinfdot}, 
\end{align}
\end{subequations}
\textcolor{change}{where $\pm a$ denotes $+a-a$.} Since the residual modes in \eqref{residue} are not explicitly accounted for in the observer-based controller \eqref{control}, the input in \eqref{residue} is treated as a disturbance. The negativity of \eqref{Vinfdot} indicates that the $L^2$ gain from the unknown input $u$ to the residue $\zeta$ is less than $\gamma$. By augmenting $\dot V_z + \dot V_e$ with the term $u^2 - \gamma^{-2} \zeta^2$, we compensate the cross terms involving the residue $\zeta$. This allows \eqref{Vzdot} and \eqref{Vedot} to be analyzed separately from the residual dynamics while avoiding spillover. 
\subsection{The $L^2$ gain for the residue}\label{sec:L2gain}
In this section, we find the smallest $\gamma$ such that \eqref{Vinfdot} is not positive. First, Lemma~\ref{lem:harmonic} gives the following harmonic bound on the residue: 
\begin{equation}\label{zetaIneq}
\textstyle\zeta^2=\frac{2}{\pi}\left(\sum_{n=N}^\infty z_n\right)^2\le\frac{2}{\pi}\sum_{n=N}^\infty\mu_nz_n^2
\end{equation}
for any $\mu_n$ satisfying the harmonic condition 
\begin{equation*}
\textstyle\sum_{n=N}^\infty\mu_n^{-1}\le1,\quad \mu_n>0. 
\end{equation*}
Using \eqref{zetaIneq}, we can find the smallest $\gamma$ such that \eqref{Vinfdot} is not positive. Namely, calculating $\dot V_\infty$ along the trajectories of \eqref{residue} and using \eqref{zetaIneq} to bound $\zeta^2$, we find 
\begin{equation*}
\begin{aligned}
\dot V_\infty+\gamma^{-2}\zeta^2{}&\textstyle-u^2+\frac{\sigma}{\gamma}\sum_{n=N}^\infty z_n^2-\frac1{\gamma\sigma}\sum_{n=N}^\infty f_n^2\\
\textstyle\le{}\sum_{n=N}^{\infty}&\textstyle\left[\frac{2}{\gamma}(q-\lambda_n)+\frac{2}{\gamma^2\pi}\mu_n+\frac{\sigma}{\gamma}\right]z_n^2\\
+&\textstyle\left[\frac{2}{\gamma}\sum_{n=N}^{\infty} z_nb_nu-u^2\right]\\
+&\textstyle\left[\frac{2}{\gamma}\sum_{n=N}^{\infty}z_nf_n-\frac{1}{\gamma\sigma}\sum_{n=N}^{\infty}f_n^2\right]. 
\end{aligned}
\end{equation*}
Completing the squares, we obtain 
\begin{multline*}
\textstyle\frac{2}{\gamma}\sum_{n=N}^{\infty} z_nb_nu-u^2\\
\textstyle=\left(\gamma^{-1}\sum_{n=N}^\infty z_nb_n\right)^2-\left(\gamma^{-1}\sum_{n=N}^\infty z_nb_n-u\right)^2\\
\textstyle\le\gamma^{-2}\left(\sum_{n=N}^\infty z_nb_n\right)^2\le\frac{2}{\gamma^2\pi}\sum_{n=N}^\infty \mu_nz_n^2. 
\end{multline*}
The last inequality follows from Lemma~\ref{lem:harmonic} since $b_n^2=2/\pi$ for $n\ge N$. Another completion of squares gives 
\begin{equation*}
\textstyle\frac{2}{\gamma}\sum_{n=N}^{\infty}z_nf_n-\frac{1}{\gamma\sigma}\sum_{n=N}^{\infty}f_n^2
\le\frac{\sigma}{\gamma}\sum_{n=N}^\infty z_n^2. 
\end{equation*}
Therefore, 
\begin{multline*}
\textstyle\dot V_\infty+\gamma^{-2}\zeta^2-u^2+\frac{\sigma}{\gamma}\sum_{n=N}^\infty z_n^2-\frac{1}{\gamma\sigma}\sum_{n=N}^\infty f_n^2\\
\textstyle\le\sum_{n=N}^{\infty}\left[\frac{4\mu_n}{\gamma^2\pi}-\frac{2}{\gamma}\left(\lambda_n-q-\sigma\right)\right]z_n^2. 
\end{multline*}
The largest $\mu_n$ such that the series is non-positive are 
\begin{equation}\label{mu_n}
    \mu_n=\gamma\frac{\pi}{2}\left(\lambda_n-q-\sigma\right),\qquad n\ge N. 
\end{equation}
Note that \eqref{N} guarantees $\mu_n>0$. The smallest $\gamma$ so that $\mu_n$ in \eqref{mu_n} satisfy the harmonic condition below \eqref{zetaIneq} is 
\begin{equation}\label{gamma:series} 
\textstyle\gamma=\frac{2}{\pi}\sum_{n=N}^{\infty}\frac{1}{\lambda_n-q-\sigma}. 
\end{equation}
The series can be calculated explicitly using the Mittag-Leffler expansion for the cotangent \cite[Section~7.10]{Spiegel2009}: 
\begin{equation*}
\textstyle\pi d\cot\pi d=1+2 \sum_{n=1}^{\infty} \frac{d^2}{d^2-n^2}. 
\end{equation*}
\textcolor{change}{For $\lambda_n=n^2$ and $d=\sqrt{q+\sigma}$, we have 
\begin{equation}\label{gamma}
\begin{aligned}
\gamma&\textstyle=\frac{1}{\pi d^2}\left[2\sum_{n=1}^{\infty}\frac{d^2}{n^2-d^2}-2\sum_{n=1}^{N-1}\frac{d^2}{n^2-d^2}\right]\\
&\textstyle=\frac{1}{\pi d^2}\left[1-\pi d\cot\pi d+2\sum_{n=1}^{N-1}\frac{d^2}{d^2-n^2}\right]. 
\end{aligned}
\end{equation}}%
Note that $d\in\N_0$ are removable singularities. Since the series \eqref{gamma:series} converges, we have $\gamma\to0$ as $N\to\infty$. This implies that as more modes are included in the control design, the destabilizing effect of the residue diminishes. 
\subsection{Spillover-free controller design}\label{sec:Stability analysis}
Using the $L^2$ residue gain calculated in the previous section, the next theorem provides a method of designing the controller and observer gains, $K$ and $L$, such that the finite-dimensional controller \eqref{control} stabilizes the semilinear heat equation \eqref{HeatBVP}. 
\begin{thm}[Stabilization of the semilinear system]\label{th:stability}
    Consider the semilinear heat equation \eqref{HeatBVP} with a continuous $f$ subject to \eqref{fSector} with some $\sigma>0$. For $N$ satisfying \eqref{N}, let $0<X\in\R^{N\times N}$ and $0<Z\in\R^{N\times N}$ be the stabilizing solutions of the algebraic Riccati equations 
    \begin{subequations}\label{ARE}
    \begin{align}
        &\textstyle XA+A^{\top}X-X\left(BB^{\top}-\gamma\sigma I_N\right)X+\frac{\sigma}{\gamma}I_N=0, \label{ARE4X}\\
        &\textstyle ZA^{\top}+AZ-Z\left(C^\top C-\gamma\sigma I_N\right)Z+\frac{\sigma}{\gamma}I_N=0 \label{ARE4Z}
    \end{align}
    \end{subequations}
    with $A$, $B$, and $C$ defined in \eqref{notations} and $\gamma$ given in \eqref{gamma}. If 
    \begin{equation}\label{spectralCondition}
        \rho(XZ)<\gamma^{-2}\quad\text{(spectral radius condition)}, 
    \end{equation}
    then the finite-dimensional dynamic controller \eqref{control} with 
    \begin{equation}\label{controlGains}
        K=B^\top X\quad\text{and}\quad L=Z(I_N-\gamma^2XZ)^{-1}C^\top
    \end{equation}
    stabilizes \textcolor{change}{the zero solution of} \eqref{HeatBVP} in the $L^2$ norm, that is, 
    \begin{equation}\label{stabilityDefinition}
    \exists M\textcolor{change}{{}\ge1}\colon\ \|z(\cdot,t)\|^2+|e^N(t)|^2\le M\|z(\cdot,0)\|^2,\ t\ge0. 
    \end{equation}
\end{thm}
{\bf Proof.} Consider the Lyapunov functional \eqref{LyapunovFunctional} and the bound on its derivative \eqref{Vdot}. In Section~\ref{sec:L2gain}, we showed that \eqref{Vinfdot} is not positive for $\gamma$ given in \eqref{gamma}. Now, we show that both \eqref{Vzdot} and \eqref{Vedot} are not positive for $X>0$ satisfying \eqref{ARE4X} and $Y=\gamma^{-2}Z^{-1}-X$ with $Z>0$ satisfying~\eqref{ARE4Z}. The spectral condition \eqref{spectralCondition} guarantees $Y>0$. Substituting \eqref{Nmodes} for $\dot z^N$ in $\dot V_z$, we find 
\begin{multline*}
\dot V_z+u^2-|B^\top Xe^N|^2+\frac{\sigma}{\gamma}|z^N|^2-\frac1{\gamma\sigma}|F|^2+\frac1{\gamma\sigma}|\tilde F|^2\\
=\left(z^N\right)^{\top}\left[XA+A^\top X+\textcolor{change}{\frac{\sigma}{\gamma}} I_N\right]z^N-|B^\top Xe^N|^2\\
+\left[2(z^N)^\top XBu+u^2\right]+\left[2(z^N)^\top XF\!-\!\frac{1}{\gamma\sigma}|F|^2\!+\!\frac1{\gamma\sigma}|\tilde F|^2\right]. 
\end{multline*}
Completing the squares, we obtain 
\begin{equation*}
\begin{aligned}
2(z^N)^\top XBu+u^2&\stackrel{\phantom{\eqref{u}}}{=}|B^{\top}Xz^N+u|^2-|B^\top Xz^N|^2\\
&\stackrel{\eqref{u}}{=}|B^{\top}Xe^N|^2-|B^\top Xz^N|^2. 
\end{aligned}
\end{equation*}
Substituting $\tilde F=F-\gamma\sigma X z^N$, we find 
\begin{equation*}
\textstyle 2\left(z^N\right)^{\top}XF-\frac{1}{\gamma\sigma}|F|^2+\frac1{\gamma\sigma}|\tilde F|^2=\gamma\sigma(z^N)^\top X^2z^N. 
\end{equation*}
Therefore, 
\begin{multline*}
\dot V_z+u^2-|B^\top Xe^N|^2+\frac{\sigma}{\gamma}|z^N|^2-\frac1{\gamma\sigma}|F|^2+\frac1{\gamma\sigma}|\tilde F|^2\\
=\left(z^N\right)^{\top}\left[XA\!+\!A^{\top}X\!-\!X(BB^\top\!-\!\gamma\sigma I_N)X\!+\!\frac{\sigma}{\gamma}I_N\right]z^N, 
\end{multline*}
which is zero by \eqref{ARE4X}. Therefore, \eqref{Vzdot} is zero. 

Consider now \eqref{Vedot}. Substituting the dynamics of $e^N(t)$ from \eqref{e} and completing the squares for $\zeta$ and $\tilde F$, we have 
\begin{equation*}
\begin{aligned}
	\dot V_e&\textstyle{}-\gamma^{-2}\zeta^2+|B^\top Xe^N|^2-\frac1{\gamma\sigma}|\tilde F|^2\\
    \stackrel{\eqref{e}}{=}{}&2\left(e^N\right)^{\top} Y\left[(A+\gamma\sigma X-LC)e^N+L\zeta-\tilde{F}\right]\\
    &\textstyle{}-\gamma^{-2}\zeta^2+\left|B^{\top}Xe^N\right|^2-\frac1{\gamma\sigma}|\tilde F|^2\\
    \le{}&\left(e^N\right)^{\top}\Bigl[Y(A+\gamma\sigma X)+(A+\gamma\sigma X)^{\top} Y-2YLC\\
    &{}+\gamma^2YLL^{\top}Y+XBB^{\top}X+\gamma\sigma Y^2\Bigr] e^{N}.
\end{aligned}
\end{equation*}
Since $L\stackrel{\eqref{controlGains}}{=}Z(I_N-\gamma^2 XZ)^{-1}C^\top=\gamma^{-2}Y^{-1}C^\top$, we have 
\begin{equation*}
    -2YLC+\gamma^2YLL^{\top}Y=-\gamma^{-2}C^\top C, 
\end{equation*}
which gives 
\begin{equation*}
\begin{aligned}
\dot V_e&\textstyle{}-\gamma^{-2}\zeta^2+|B^\top Xe^N|^2-\frac1{\gamma\sigma}|\tilde F|^2\\
\le{}&\left(e^N\right)^{\top}\Bigl[Y(A+\gamma\sigma X)+(A+\gamma\sigma X)^{\top} Y\\
&{}+XBB^{\top}X+\gamma\sigma Y^2-\gamma^{-2}C^\top C\Bigr] e^{N}\\
=&\left(e^N\right)^{\top}[\gamma^{-2}Z^{-1}\times\eqref{ARE4Z}\times Z^{-1}-\eqref{ARE4X}]e^N. 
\end{aligned}
\end{equation*}
The latter is zero by \eqref{ARE}. Therefore, $\dot V(t)\le0$, which implies $V(t)\le V(0)$. Parseval's theorem and \eqref{LyapunovFunctional} give  
\begin{equation*}
\begin{aligned}
    &c_1\left(\|z(\cdot,t)\|^2\!+\!|e^N(t)|^2\right)\le V(t)\le c_2\left(\|z(\cdot,t)\|^2\!+\!|e^N(t)|^2\right),\\
    &0<c_1=\min\{\lambda_{\min}(X),\lambda_{\min}(Y),\gamma^{-1}\},\\
    &0<c_2=\max\{\lambda_{\max}(X),\lambda_{\max}(Y),\gamma^{-1}\}, 
\end{aligned}
\end{equation*}
where $\lambda_{\min}$ and $\lambda_{\max}$ are the minimum and maximum eigenvalues, respectively. Therefore, 
\begin{multline*}
    \|z(\cdot,t)\|^2+|e^N(t)|^2\le c_1^{-1}V(t)\\
    \le c_1^{-1}V(0)\le c_1^{-1}c_2\left(\|z(\cdot,0)\|^2+|e^N(0)|^2\right). 
\end{multline*}
Since $\hat z^N(0)=0$, we have $|e^N(0)|^2=|z^N(0)|^2$, which gives \eqref{stabilityDefinition} with $M=2c_1^{-1}c_2$. \hfill$\qed$
\begin{rem}[Exponential stability]\label{rem:expStabl}
The modified state $z_\alpha(x,t)=e^{\alpha t}z(x,t)$ satisfies \eqref{HeatBVP} with $q$ replaced by $q+\alpha$ and $f$ replaced by $e^{\alpha t}f(\cdot,t,e^{-\alpha t}z_\alpha(\cdot,t))$, which satisfies \eqref{fSector} with the same Lipschitz constant $\sigma$. Therefore, the exponential stability of \eqref{HeatBVP} in the $L^2$ norm with the decay rate $\alpha$ follows from Theorem~\ref{th:stability} if $q$ is replaced by $q+\alpha$. 
\end{rem}{\color{change}
\begin{prop}[Feasibility for small $\sigma$]\label{prop:sigma}
Let $q\neq\lambda_n$ for all $n \in \N_0$. Then, for a large enough $N$ and small enough $\sigma>0$, \eqref{ARE} always admit stabilizing solutions $X>0$ and $Z>0$ that satisfy \eqref{spectralCondition}.
\end{prop}}
\textcolor{change}{The proof is given in Appendix~\ref{prop:sigma:proof}. Note that $q \neq \lambda_n$ can always be guaranteed for all $n \in \N_0$ by considering exponential stability with a sufficiently small decay rate~$\alpha>0$, as explained in Remark~\ref{rem:expStabl}.}
\textcolor{change}{\begin{rem}[Linear case]\label{cor:linear}
In the proof of Proposition~\ref{prop:sigma}, we showed that \eqref{ARE} and \eqref{spectralCondition} with $\sigma=0$ have the solution
\[
    X=\operatorname{diag}\{X_0,0_{N-N_0}\}, 
    \quad 
    Z=\operatorname{diag}\{Z_0,0_{N-N_0}\},
\]
where $X_0>0$ and $Z_0>0$ are of size $N_0$, chosen so that $\lambda_{N_0-1} < q < \lambda_{N_0}$, and $N$ is large enough to ensure $\rho(XZ) < \gamma^{-2}$. Substituting into \eqref{controlGains} yields
\[
K = \bigl[B_0^\top X_0,\; 0_{1 \times (N - N_0)} \bigr],
\quad
L=\left[\begin{smallmatrix}
Z_0(I_{N_0}-\gamma^2X_0Z_0)^{-1}C_0^\top\\ 
0_{(N-N_0)\times 1}
\end{smallmatrix}\right], 
\]
where $B_0$ is defined in \eqref{notations} with $N$ replaced by $N_0$. Thus, Theorem~\ref{th:stability} indicates that in the linear case, the controller and observer gains should be set to zero from index $N_0$ onward, where $N_0$ corresponds to the number of unstable modes. This provides a formal justification for the common practice of zeroing certain gains, as discussed in \cite{Katz2020a,Katz2021b,Lhachemi2022a,Wang2023}.
\end{rem}}


\textcolor{change}{If $\sigma \neq 0$ is not small, feasibility is not guaranteed, even for large $N$. Intuitively, increasing $N$ might appear to restrict the admissible range of $\sigma$, as it amplifies the effect of the nonlinearity $F$ in \eqref{Nmodes}, potentially leading to spillover. However, the next proposition proves that this is not the case: increasing the controller order $N$ can only expand the set of admissible $\sigma$ values.}

\begin{prop}[Feasibility for larger $N$]\label{prop:N}
If the conditions of Theorem~\ref{th:stability} hold for some $N$, then they hold for $N+1$. In particular, the maximum Lipschitz constant, $\sigma$, such that \textcolor{change}{the zero solution of \eqref{HeatBVP} is stable under} \eqref{control} does not decrease when the number of modes $N$ used by the controller increases. 
\end{prop}
\textcolor{change}{The proof is given in Appendix~\ref{prop:N:proof}.}%
\section{Sample-and-hold control of the semilinear heat equation via the $L^2$ residue separation}\label{sec:sampled-data}
\textcolor{change}{An output-feedback controller for \eqref{HeatBVP} with known nonlinearity was proposed in \cite{Katz2022b} using the lifting transformation \( w(x,t) = z(x,t) - \psi(x)u(t) \), where \(\psi\) is the eigenfunction of \(\partial_x^2\) satisfying \(\psi'(0) = 0\) and \(\psi'(\pi) = 1\). While this transformation homogenizes the boundary conditions, it introduces the time derivative of the input into the PDE. In contrast, the \(L^2\)-residue separation method enabled by the harmonic inequality (Lemma~\ref{lem:harmonic}) eliminates the need for lifting. This not only reduces the controller’s order but also enables sampled-data implementation using a standard zero-order hold, rather than the generalized hold required in, e.g., \cite{Katz2021a}. In this section, we present the design of such a sampled-data controller.}


Consider the semilinear heat equation with sample-and-hold input 
\begin{subequations}\label{HeatBVPdelay}
\begin{align}
& z_t=z_{xx}+qz+f(\cdot,t,z(\cdot,t)),\label{PDEdelay}\\
& z_x(0,t)=0, \quad z_x(\pi, t)=u(t_k),\quad t\in[t_k,t_{k+1}),\label{BCdelay}\\
& y(t)=z(0,t), \label{Outputdelay}
\end{align}
\end{subequations}
with sampling instants $\{t_k\}_{k\in\N_0}$ satisfying 
\begin{equation*}
    0=t_0<t_1<t_2<\cdots,\quad \sup_{k\in\N_0}(t_{k+1}-t_k)\le h. 
\end{equation*}
Performing modal decomposition for \eqref{HeatBVPdelay} similarly to Section~\ref{sec:modalDecomposition}, we obtain (cf.~\eqref{HeatModalDecomposition})
\begin{subequations}\label{HeatModalDecompositionDelay}
\begin{align}
\dot{z}^N(t)&=Az^N(t)+Bu(t_k)+F(t),\label{NmodesDelay}\\
\dot{z}_n(t)&=-\left(\lambda_n-q\right) z_n+b_nu(t_k)+f_n(t),\ n\ge N, \label{residueDelay}\\
y(t)&=Cz^N(t)+\zeta(t) \label{outputDelay}
\end{align}
\end{subequations}
when $t\in[t_k,t_{k+1})$, with the notations from \eqref{notations}. Then, the observer-based controller \eqref{control} should be changed to 
\begin{subequations}\label{controlDelay}
\begin{align}
\dot{\hat z}^N(t)&=(A+\gamma\sigma X)\hat z^N(t)+Bu(t_k)\nonumber\\
&\hspace{1cm}{}-L(y(t)-C\hat z^N(t)),\quad t\in[t_k,t_{k+1}),\label{observerDelay}\\
u(t_k)&=-K\hat z^N(t_k)\label{uDelay}
\end{align}
\end{subequations}
with $\hat z^N(0)=0$. The dynamics of the estimation error $e^N(t)=\hat z^N(t)-z^N(t)$ is governed by the same equation \eqref{e} because the sample-and-hold input cancels out. 

The well-posedness of \eqref{HeatBVPdelay}, \eqref{controlDelay} is established using the step method. If $z(\cdot,0)\in H^1(0,\pi)$, by Theorems~6.1.2 and 6.3.1 of \cite{Pazy1983}, a unique classical solution exists on $[0,t_1)$: 
\begin{equation*}
\begin{aligned}
    &z\in C([0,t_1),L^2)\cap C^1((0,t_1),L^2),\\
    &z(\cdot,t)\in H^2(0,\pi),\quad\forall t\in(0,t_1]. 
\end{aligned}
\end{equation*}
Since $z(\cdot,t_1)\in H^2(0,\pi)$, the same theorems guarantee that \eqref{HeatBVPdelay}, \eqref{controlDelay} has a unique classical solution on $[t_1,t_2)$ for a constant $u(t_1)$. Using the same reasoning sequentially on each sampling interval $[t_k,t_{k+1})$, we establish the existence of the unique classical solution 
\begin{equation*}
\begin{aligned}
    &z\in C([0,\infty),L^2)\cap C^1((0,\infty)\setminus\mathcal{J},L^2),\quad\mathcal{J}=\{t_k\}_{k\in\N},\\
    &z(\cdot,t)\in H^2(0,\pi),\quad\forall t>0. 
\end{aligned}
\end{equation*}

To derive the stability conditions, we ensure that the observer state $\hat z$, used by the controller in \eqref{uDelay}, does not change too rapidly between sampling instants, that is, $u(t_k)\approx -K\hat z^N(t)$ on $[t_k,t_{k+1})$. To achieve this, we need to examine the dynamics of $\hat z^N = z^N + e^N$, which couples the subsystem states $z^N$ and $e^N$. This coupling makes the separation \textcolor{change}{in the Lyapunov analysis}, as seen in \eqref{Vzdot} and \eqref{Vedot}, ineffective. Instead, we construct a Lyapunov--Krasovskii functional and analyze it for both subsystems simultaneously, which leads to stability conditions expressed in terms of linear matrix inequalities (LMIs) rather than algebraic Riccati equations.
\begin{figure*}
\begin{equation}\label{Psi}
\begin{aligned}
&\Psi=\left[\begin{smallmatrix}
\Psi_{11} & \Psi_{12} & \gamma\sigma P_z & 0 & K^\top K-P_zBK & K^\top K-P_zBK & h(A-BK)^\top W_z & h\gamma\sigma XW_e\\
* & \Psi_{22} & -\gamma\sigma P_e & P_eL & K^\top K & K^\top K & -hK^\top B^\top W_z & h(A-LC+\gamma\sigma X)^\top W_e \\
* & * & -\gamma\sigma I_N & 0 & 0 & 0 & h\gamma\sigma W_z & -h\gamma\sigma W_e \\
* & * & * & -\gamma^{-2} & 0 & 0 & 0 & hL^\top W_e \\
* & * & * & * & -\frac{\pi^2}4W_z+K^\top K & K^\top K & -hK^\top B^\top W_z & 0 \\
* & * & * & * & * & -\frac{\pi^2}4W_e+K^\top K & -hK^\top B^\top W_z & 0 \\
* & * & * & * & * & * & -W_z & 0 \\
* & * & * & * & * & * &  * & -W_e
\end{smallmatrix}\right]\\
&\Psi_{11}=P_z(A-BK)+(A-BK)^\top P_z+K^\top K+\frac{\sigma}{\gamma}I_N,\\
&\Psi_{12}=-P_zBK+\gamma\sigma XP_e+K^\top K,\\
&\Psi_{22}=P_e(A-LC+\gamma\sigma X)+(A-LC+\gamma\sigma X)^\top P_e+K^\top K. 
\end{aligned}
\end{equation}
\end{figure*}
\begin{thm}[Stability under sample-and-hold]\label{th:inputDelay}
Consider the semilinear heat equation \eqref{HeatBVPdelay} with a continuous $f$ subject to \eqref{fSector} with some $\sigma>0$. For $N$ satisfying \eqref{N}, let $0<X\in\R^{N\times N}$ and $0<Z\in\R^{N\times N}$ be the stabilizing solutions of \eqref{ARE}, \eqref{spectralCondition} with $A$, $B$, and $C$ defined in \eqref{notations} and $\gamma$ given in \eqref{gamma}. Let there exist positive-definite $N\times N$ matrices $P_z$, $P_e$, $W_z$, and $W_e$, such that $\Psi\le0$ with $\Psi$ defined in \eqref{Psi}. Then the dynamic controller \eqref{controlDelay} with $K$ and $L$ given by \eqref{controlGains} stabilizes \textcolor{change}{the zero solution of} \eqref{HeatBVPdelay} in the $L^2$ norm in the sense of \eqref{stabilityDefinition}. 
\end{thm}
{\bf Proof.} 
Define the sampling-induced errors 
\begin{equation*}
    \delta_z(t)=z^N(t_k)-z^N(t)\quad\text{and}\quad\delta_e(t)=e^N(t_k)-e^N(t)
\end{equation*}
for $t\in[t_k,t_{k+1})$. Then 
\begin{equation}\label{controlSDCont}
\begin{aligned}
u(t_k)&=-K\hat z^N(t_k)=-K[z^N(t_k)+e^N(t_k)]\\
&=-K[z^N(t)+\delta_z(t)+e^N(t)+\delta_e(t)], 
\end{aligned}
\end{equation}
and the closed-loop system \eqref{NmodesDelay}, \eqref{uDelay} takes the form  
\begin{multline}\label{PlantSD-CL}
    \dot z^N(t)=(A-BK)z^N(t)-BKe^N(t)+F(t)\\
    -BK\delta_z(t)-BK\delta_e(t). 
\end{multline}
Consider the Lyapunov functional 
\begin{equation}\label{LyapunovFunctionalDelay}
V_h=V_{P,z}+V_{P,e}+V_{W,z}+V_{W,e}+V_\infty,
\end{equation}
where 
\begin{equation*}
\begin{aligned}
&V_{P,z}=|z^N(t)|_{P_z}^2, \\
&V_{P,e}=|e^N(t)|_{P_e}^2, \\
&\textstyle V_{W,z}=h^2\int_{t_k}^t|\dot z^N(s)|_{W_z}^2ds-\frac{\pi^2}4\int_{t_k}^t|\delta_z(s)|_{W_z}^2ds,\\
&\textstyle V_{W,e}=h^2\int_{t_k}^t|\dot e^N(s)|_{W_e}^2ds-\frac{\pi^2}4\int_{t_k}^t|\delta_e(s)|_{W_e}^2ds,\\
\end{aligned}
\end{equation*}
and $V_\infty$ is as in \eqref{LyapunovFunctional}. The terms $V_{W,z}$ and $V_{W,e}$, inspired by \cite{Liu2012}, \textcolor{change}{do not grow at $t_k$ and are nonnegative. They }will bound the sampling induced errors $\delta_z$ and $\delta_e$ using the dynamics of $z^N$ and $e^N$. 

For $t\in[t_k,t_{k+1})$, we have  (cf.~\eqref{Vdot})
\begin{subequations}\label{VdotDelay}
\begin{align}
&\dot V_h\le\dot V_{P,z}+\dot V_{P,e}+\dot V_{W,z}+\dot V_{W,e}+\dot V_\infty\notag\\
&\hspace{3.8cm}\pm u^2(t_k)\pm\gamma^{-2}\zeta^2+\eqref{SProc_f}\notag\\
&=\Bigl[\dot V_{P,z}+\dot V_{P,e}+\dot V_{W,z}+\dot V_{W,e}\notag\\
&\hspace{1.65cm}{}+u^2(t_k)-\gamma^{-2}\zeta^2+\frac{\sigma}{\gamma}|z^N|^2-\frac1{\gamma\sigma}|F|^2\Bigr]\label{Vzedot}\\
&+\Bigl[\dot V_\infty-\!u^2(t_k)\!+\!\gamma^{-2}\zeta^2\!+\frac{\sigma}{\gamma}\sum_{n=N}^\infty z_n^2-\frac{1}{\gamma\sigma}\sum_{n=N}^\infty f_n^2\Bigr]\label{VinfdotDelay}. 
\end{align}
\end{subequations}
The $L^2$ gain calculation of Section~\ref{sec:L2gain} remains valid with $u(t)$ replaced by $u(t_k)$, implying that \eqref{VinfdotDelay} is not positive for the same $\gamma$ as in \eqref{gamma:series}. Below, we show that \eqref{Vzedot} is also not positive.  

We have  
\begin{multline*}
\dot V_{P,z}\stackrel{\eqref{PlantSD-CL}}{=}2(z^N)^\top P_z[(A-BK)z^N-BKe^N+F\\
    -BK\delta_z-BK\delta_e],
\end{multline*}
\begin{multline*}
\dot V_{P,e}\stackrel{\eqref{e}}{=}2(e^N)^\top P_e[(A-LC+\gamma\sigma X)e^N\\
+\gamma\sigma Xz^N-F+L\zeta], 
\end{multline*}
\begin{equation*}
\begin{aligned}
&\textstyle\dot{V}_{W,z}=h^2\left(\dot{z}^N(t)\right)^{\top} W_z \dot{z}^N(t)-\frac{\pi^2}{4} \delta_z^{\top}(t) W_z \delta_z(t), \\
&\textstyle\dot{V}_{W,e}=h^2\left(\dot{e}^N(t)\right)^{\top} W_e \dot{e}^N(t)-\frac{\pi^2}{4} \delta_e^{\top}(t) W_e \delta_e(t).
\end{aligned}
\end{equation*}
Substituting the above and \eqref{controlSDCont} into \eqref{Vzedot}, we obtain 
Substituting \eqref{controlSDCont} into \eqref{Vzedot}, we obtain 
\begin{equation*}
\begin{aligned}
&\dot V_{P,z}+\dot V_{P,e}+\dot V_{W,z}+\dot V_{W,e}\\
&\textstyle\hspace{2.8cm}{}+u^2(t_k)-\gamma^{-2}\zeta^2+\frac{\sigma}{\gamma}|z^N|^2-\frac1{\gamma\sigma}|F|^2\\
&\le\psi^\top \tilde\Psi\psi+h^2\left(\dot{z}^N\right)^{\top} W_z \dot{z}^N+h^2\left(\dot{e}^N\right)^{\top} W_e \dot{e}^N, 
\end{aligned}    
\end{equation*}
where $\psi=\operatorname{col}\{z^N,e^N,(\gamma\sigma)^{-1}F,\zeta,\delta_z,\delta_e\}$ and $\tilde\Psi$ is obtained by removing the last two block rows and columns from $\Psi$. Substituting \eqref{PlantSD-CL} for $\dot z$ and \eqref{e} for $\dot e$, and using the Schur complement lemma, we find that $\Psi\le0$ guarantees that \eqref{Vzedot} is not positive. Therefore, $\dot V_h\le0$. Though $V_h(t)$ may be discontinuous at $t_k$, it does not grow at these points since $V_{W,z}(t_k)=0=V_{W,e}(t_k)$. The remainder of the proof is similar to that of Theorem~\ref{th:stability}. 
\hfill$\qed$

Theorem~\ref{th:inputDelay} deals with the sample-and-hold implementation of the controller from Theorem~\ref{th:stability}. The following result shows that in the non-critical case, when \eqref{ARE} remains feasible for a slightly larger $q$, the feasibility of Theorem~\ref{th:inputDelay} for small ${h>0}$ follows from the feasibility of Theorem~\ref{th:stability}. In other words, if the continuous-time controller stabilizes the system, it will also stabilize it under sufficiently fast sampling.
\begin{prop}[Feasibility for small $h$]\label{prop:h}
Let the conditions of Theorem~\ref{th:stability} be true for some $q>0$. If $q$ is replaced with $\tilde q\in[0,q)$ in $A$ that appears in \eqref{Psi}, then $\Psi\le0$ is feasible for a small enough sampling period $h>0$. 
\end{prop}
The proof is given in Appendix~\ref{prop:h:proof}. 
\section{Numerical examples}
\subsection{$L^2$ gain reduction: Harmonic vs Sobolev's inequalities}
\begin{figure}\centering
\includegraphics[width=\linewidth]{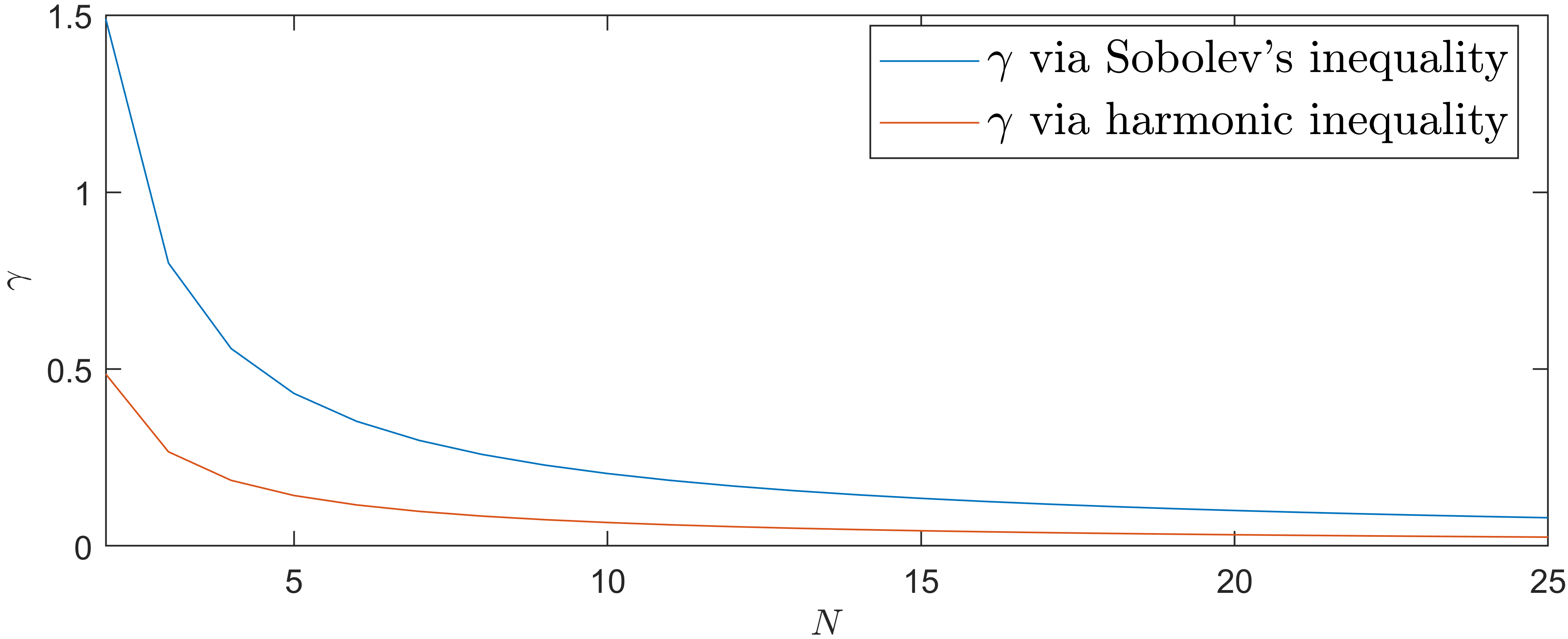}
\caption{The $L^2$ residue gain for $q=1.1$, $\sigma=0$, and $N=2,\ldots,25$ calculated using Sobolev's inequality (blue) and harmonic inequality (orange).}\label{fig:gammas}
\end{figure}
\textcolor{change}{
We employed the harmonic inequality (Lemma~\ref{lem:harmonic}) to derive~\eqref{zetaIneq}, which enabled the computation of the $L^2$ residue gain \eqref{gamma:series} by minimizing $\gamma$ subject to the constraints
\begin{equation}\label{minProblem}
\textstyle\mu_n \le \gamma \frac{\pi}{2}(\lambda_n - q - \sigma),\quad \forall n \ge N, \qquad \sum_{n=N}^\infty \mu_n^{-1} \le 1.
\end{equation}
Alternatively, Sobolev's inequality \cite[Lemma~4.1]{Kang2019} can be used to obtain (cf.~(2.33) in \cite{Katz2022b}):
\begin{equation*}
\textstyle\zeta^2(t) \le \sum_{n=N}^\infty \kappa_n z_n^2(t), \quad \kappa_n = \frac{1}{\pi} + \Gamma + \Gamma^{-1} \lambda_n, \quad \Gamma > 0.
\end{equation*}
This is a special case of \eqref{zetaIneq} with $\mu_n = \frac{\pi}{2} \kappa_n$. Minimizing $\gamma$ under \eqref{minProblem} with $\mu_n = \frac{\pi}{2} \left( \frac{1}{\pi} + \Gamma + \Gamma^{-1} \lambda_n \right)$ yields
\begin{equation}\label{gamma:Sobolev}
\textstyle\gamma = \frac{2\sqrt{\lambda_N} + \frac{1}{\pi}}{\lambda_N - q - \sigma} = \frac{2N + \frac{1}{\pi}}{N^2 - q - \sigma}, \qquad \Gamma = \sqrt{\lambda_N}.
\end{equation}
The resulting gain is larger than the one obtained via the harmonic inequality due to the additional constraint $\mu_n = \frac{\pi}{2} \left( \frac{1}{\pi} + \Gamma + \Gamma^{-1} \lambda_n \right)$. In particular, }
Figure~\ref{fig:gammas} compares the $L^2$ gains obtained using Sobolev's and harmonic inequalities for $q=1.1$, $\sigma=0$, and $N=2,\ldots,25$ (note that $N=1$ violates \eqref{N}). Both gains approach zero as $N$ increases, indicating that the $L^2$ residue gain decreases when more modes are included in the controller design. However, the harmonic inequality reduces the $L^2$ gain by a factor of 3, significantly lowering the number of modes required to analytically guarantee spillover avoidance. Specifically, \eqref{ARE} with $\gamma$ as defined in \eqref{gamma:Sobolev} requires $N = 20$ modes to be feasible, whereas $\gamma$ from \eqref{gamma} requires only $N = 7$ modes. Thus, the harmonic inequality introduced in Lemma~\ref{lem:harmonic} enhances the $L^2$ residue separation approach, reducing the number of modes needed to guarantee spillover avoidance analytically. 
\subsection{\color{change}Increasing controller order enlarges the admissible class of nonlinearities}
\begin{figure}\centering
\includegraphics[width=.8\linewidth]{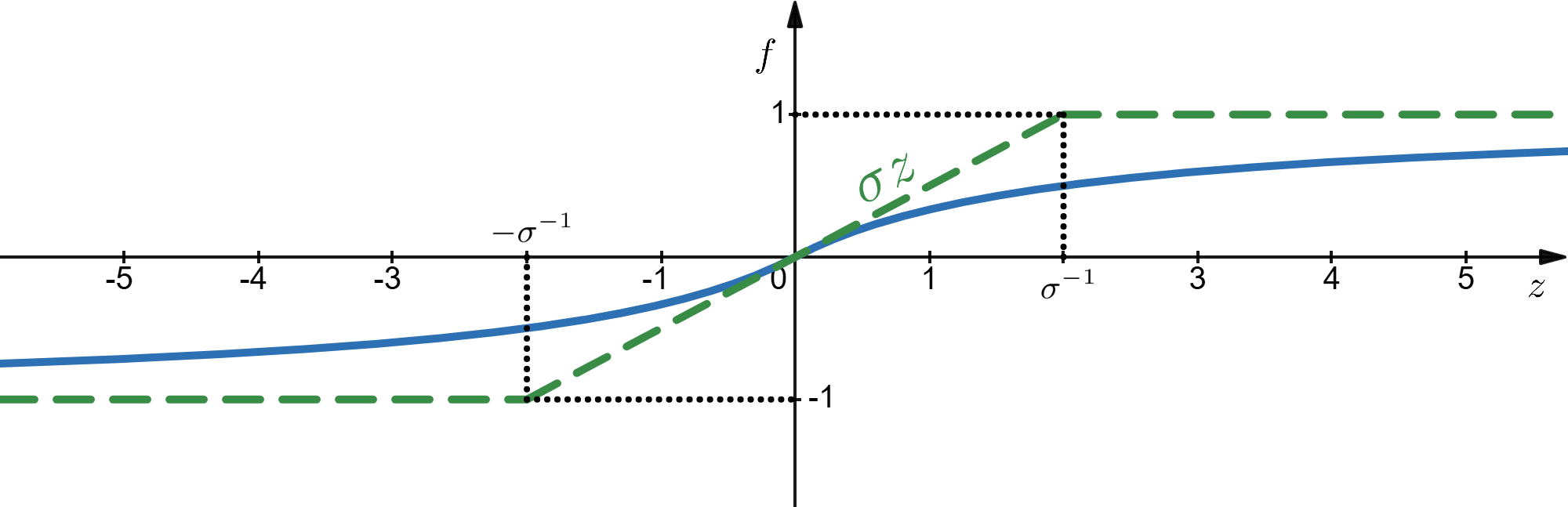}
\caption{The saturation of the linear reaction rate $\sigma z$ with the saturation range $[-\sigma^{-1},\sigma^{-1}]$ (green dashed line) and its smooth approximation given by \eqref{saturation} (blue solid line) for $\sigma=0.5$}\label{fig:saturation}
\end{figure}

Consider the semilinear heat equation \eqref{HeatBVP} with the saturated reaction rate (Fig.~\ref{fig:saturation}) 
\begin{equation}\label{saturation}
\textstyle f(x,t,z)=\frac{\sigma z}{1+\sigma|z|}, 
\end{equation}
which satisfies the sector condition \eqref{fSector}. Without control input, the system is unstable for any $q>0$ and $\sigma\ge0$. Let $q = 0.1$. Using binary search, we determine the maximum value of $\sigma$ for which the conditions of Theorem~\ref{th:stability} are feasible. The corresponding values of $\sigma$ and $\gamma$ (the $L^2$ residue gain) for various numbers of modes $N$ are listed in the table below.

\begin{small}\begin{tabular}{c|cccccc}
$N$ & $1$ & $2$ & $3$ & $4$ & $5$ & $6$ \\ \hline
$\sigma$ & $0.037$ & $0.193$ & $0.278$ & $0.327$ & $0.360$ & $0.382$ \\ 
$\gamma$ & $1.156$ & $0.427$ & $0.256$ & $0.183$ & $0.142$ & $0.116$ 
\end{tabular}
\end{small}

Observe that the value of $\sigma$ preserving stability increases as the number of controlled modes increases, consistent with Proposition~\ref{prop:N}. Meanwhile, the corresponding $L^2$ residue gains $\gamma$ decrease, indicating that the residual state becomes less influenced by the input as more modes are incorporated into the controller design.
\subsection{\color{change}Spillover-free control of the semilinear heat equation}
We now demonstrate that the sample-and-hold implementation of the control signal preserves stability, provided the sampling period is sufficiently small. Consider the nonlinear heat equation \eqref{HeatBVPdelay} with the saturated reaction rate \eqref{saturation} and a sample-and-hold control input. Theorem~\ref{th:stability} is feasible for $N=3$, $q=0.1$, and $\sigma=0.2$. Thus, the continuous-time input stabilizes the system with control gains:
\begin{equation*}
    K\approx[1.33, -0.16, 0.06]\quad\text{and}\quad L\approx[2.82, 0.01, 0.05]^\top, 
\end{equation*}
given by \eqref{controlGains}. Furthermore, Proposition~\ref{prop:h} ensures that the LMI of Theorem~\ref{th:inputDelay} is feasible for the reduced $q=0.08$ and a sufficiently small sampling period $h$. Using binary search, we find the maximum allowable sampling period $h\approx 0.1$. 

\begin{figure}\centering
\includegraphics[width=\linewidth]{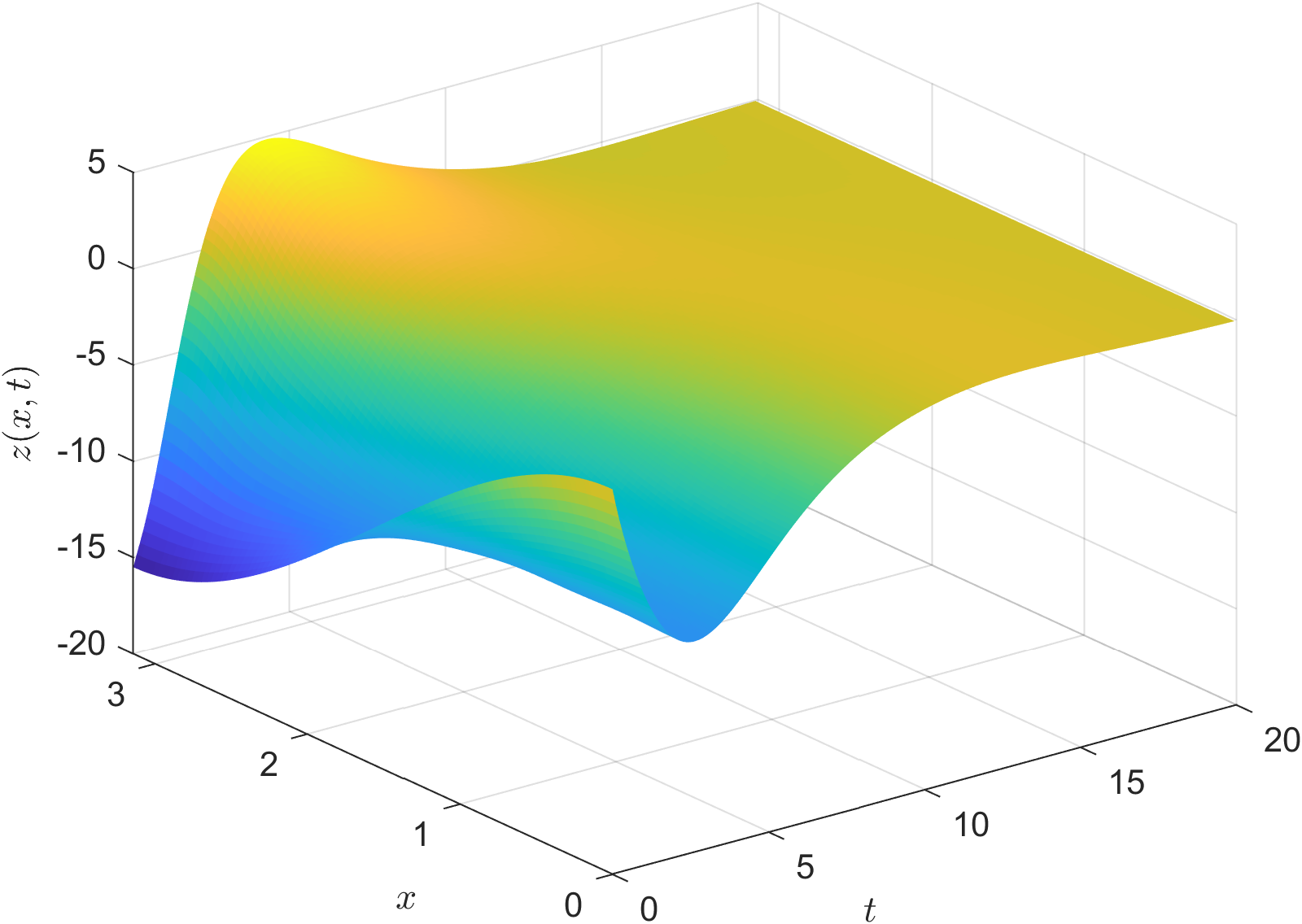}
\caption{The state $z(x,t)$ under sample-and-hold control}\label{fig:state}
\end{figure}

\begin{figure}\centering
\includegraphics[width=\linewidth]{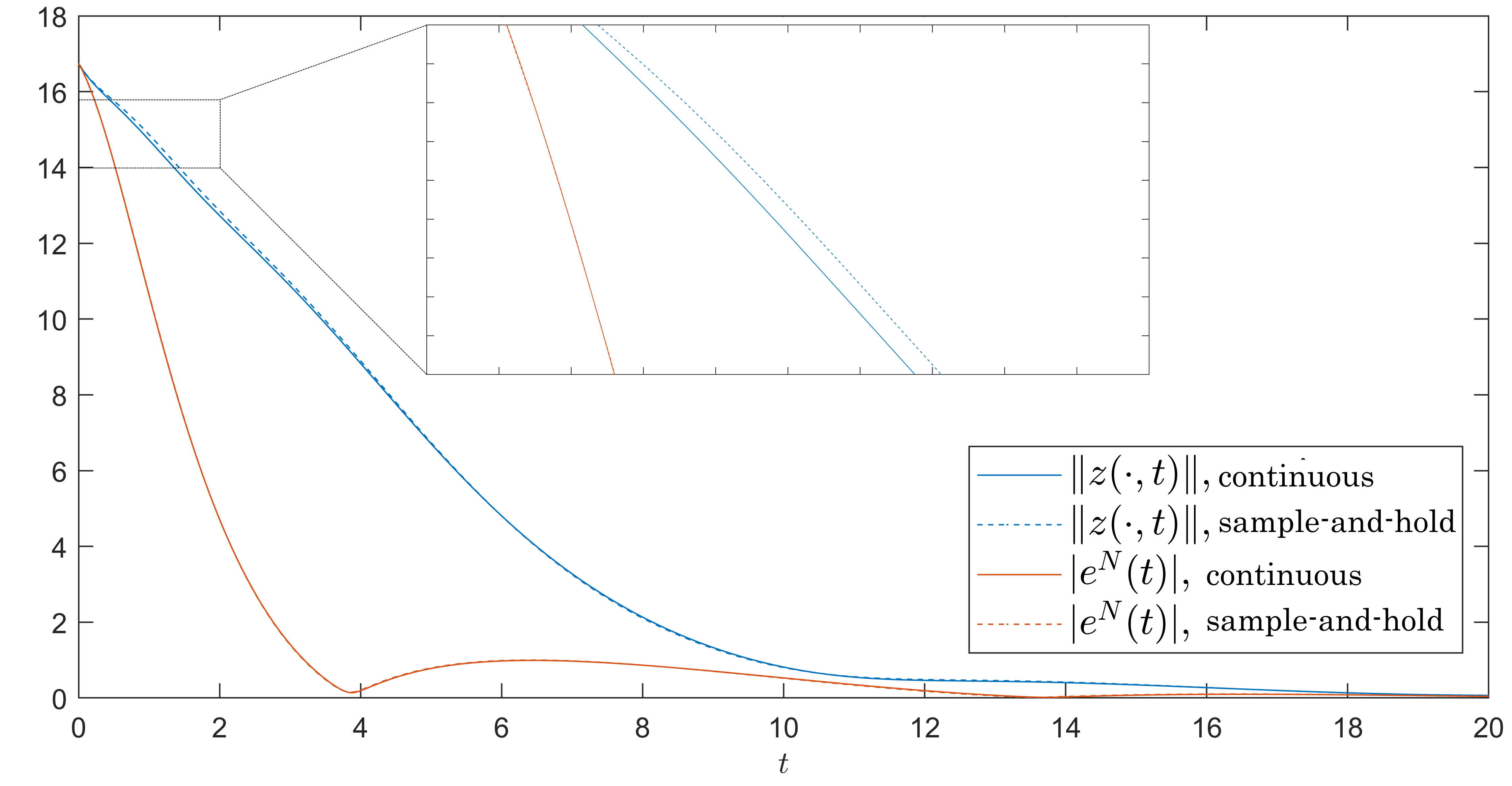}
\caption{The norms of the state (blue) and observation error (orange) for continuous (solid line) and sample-and-hold (dashed line) control}\label{fig:norms}
\end{figure}

\begin{figure}\centering
\includegraphics[width=\linewidth]{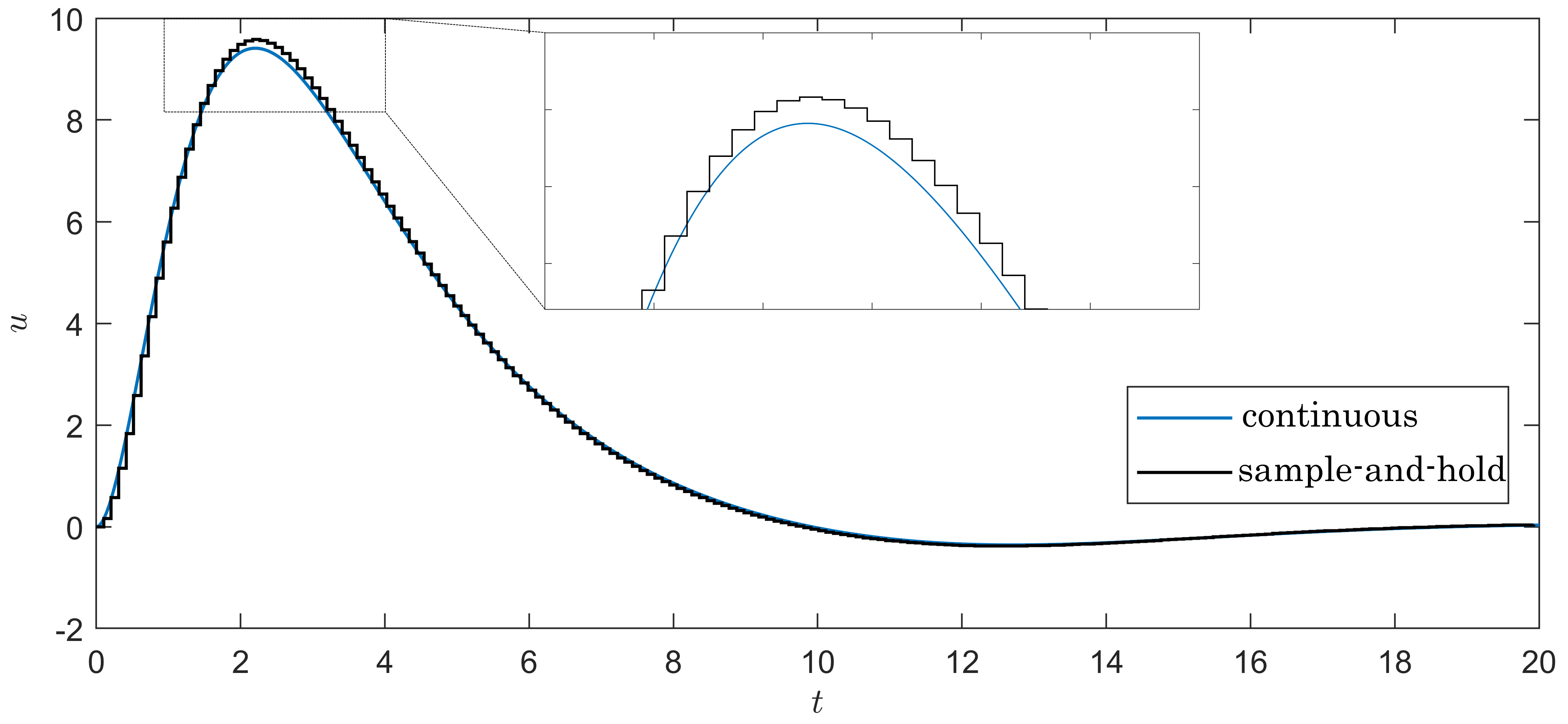}
\caption{Continuous (blue) and sample-and-hold (black) control}\label{fig:control}
\end{figure}

Numerical simulation results are presented in Figures~\ref{fig:state}--\ref{fig:control}. The initial condition was 
\begin{equation}\label{IC}
    \textstyle z(x,0)=x^3-\frac{3\pi}2x^2, \quad x\in[0,\pi],
\end{equation}
which was selected to satisfy the boundary conditions $z_x(0,0)=z_x(\pi,0)=0$. Both continuous and sample-and-hold control strategies ensure that the state and observation error converge to zero, with the control input also approaching zero.

\section{Conclusion}
We have developed the $L^2$ residue separation method for designing spillover-free \textcolor{change}{finite-dimensional} output-feedback controllers. \textcolor{change}{The newly introduced harmonic inequality enables efficient separation of the infinite-dimensional residue present in the measured output, enhancing the overall analysis. The $L^2$ residue separation framework offers several advantages: it accommodates unknown nonlinearities, provides systematic design of both the controller and observer gains, and eliminates the need for a lifting transformation, thereby supporting sampled-data control via zero-order hold. Additionally, it offers the first theoretical justification for the widely used controller and observer gain structure in the linear case. These results position $L^2$ residue separation as a powerful and versatile tool for output-feedback design in infinite-dimensional systems. Extending this approach to broader and more general classes of systems represents a promising direction for future research.}
\providecommand{\noopsort}[1]{}

\appendix
\color{change}

\section{Proof of Proposition~\ref{prop:sigma}}\label{prop:sigma:proof}
Take $N_0$ such that $\lambda_{N_0-1} < q < \lambda_{N_0}$, and define $A_0$, $B_0$, and $C_0$ as in~\eqref{notations} with $N$ replaced by $N_0$. Since $A_0$ is anti-Hurwitz, $(A_0,B_0)$ is controllable, and $(A_0,C_0)$ is observable, the matrices
\begin{equation}\label{X0Z0}
\begin{aligned}
    \textstyle X_0&=\left[\int_0^\infty e^{-A_0t}B_0B_0^\top e^{-A_0^\top t}\,dt\right]^{-1}, \\
    \textstyle Z_0&=\left[\int_0^\infty e^{-A_0^\top t}C_0^\top C_0e^{-A_0t}\,dt\right]^{-1},
\end{aligned}
\end{equation}
are well-defined, positive definite, and satisfy
\begin{equation*}
\begin{aligned}
& X_0 A_0 + A_0^\top X_0 - X_0 B_0 B_0^\top X_0 = 0,\\
& Z_0 A_0^\top + A_0 Z_0 - Z_0 C_0^\top C_0 Z_0 = 0. 
\end{aligned}
\end{equation*}
Then 
\begin{equation*}
    X=\operatorname{diag}\{X_0,0_{N-N_0}\}\quad\text{and}\quad Z=\operatorname{diag}\{Z_0,0_{N-N_0}\}
\end{equation*}
satisfy~\eqref{ARE} with $\sigma=0$. Moreover, $\rho(XZ) = \rho(X_0 Z_0) < \gamma^{-2}$ for a large enough $N\ge N_0$, where $\gamma$ defined in \eqref{gamma} goes to $0$ when $N\to\infty$. 
By \cite[Theorem~9.1.1]{Lancaster1995}, there exist $\tilde X>0$ and $\tilde Z>0$ satisfying 
\begin{equation*}
\begin{aligned}
&\tilde XA+A^\top \tilde X-\tilde XBB^\top \tilde X=-\varepsilon I<0,\\
&\tilde ZA^\top+A\tilde Z-\tilde ZC^\top C\tilde Z=-\varepsilon I<0,  
\end{aligned}
\end{equation*}
with any $\varepsilon>0$. Since $\tilde X$ and $\tilde Z$ continuously depend on $\varepsilon>0$ \cite[Theorem~2.1]{Ran1988}, \eqref{spectralCondition} implies $\rho(\tilde X\tilde Z)<\gamma^{-2}$ for a small enough $\varepsilon>0$. Then, for a small enough $\sigma$, 
\begin{equation}\label{ARE:tilde}
\begin{aligned}
&\tilde XA+A^\top \tilde X-\tilde XBB^\top \tilde X+\sigma(\gamma \tilde X^2+\gamma^{-1}I_N)<0,\\
&\tilde ZA^\top+A\tilde Z-\tilde ZC^\top C\tilde Z+\sigma(\gamma \tilde Z^2+\gamma^{-1}I_N)<0, 
\end{aligned}
\end{equation}
which are \eqref{ARE} with ``$=$'' replaced by ``$<$''. 

Now we show that \eqref{ARE:tilde} imply the feasibility of \eqref{ARE} together with \eqref{spectralCondition}. To this end, we demonstrate that \eqref{ARE:tilde} ensure that the $L^2$ gain of a suitably constructed system does not exceed $\gamma^{-1}$. We then invoke \cite[Theorem~8.3.2]{Green2012}, which guarantees the feasibility of the corresponding algebraic Riccati equations with the spectral condition. 

Let $\tilde Y=\gamma^{-2}\tilde Z^{-1}-\tilde X$, which is positive definite for $\gamma$ satisfying $\rho(\tilde X\tilde Z)<\gamma^{-2}$. Consider 
\begin{equation*}
    \tilde V=|z^N|_{\tilde X}^2+|e^N|_{\tilde Y}^2. 
\end{equation*}
Then \eqref{ARE:tilde} guarantee that \eqref{Vzdot}, \eqref{Vedot} with $X$ and $Y$ replaced by $\tilde X$ and $\tilde Y$, are not positive. Therefore, 
\begin{equation*}
    \dot{\tilde V}\le u^2+\frac{\sigma}{\gamma}|z^N|^2-\frac1{\gamma\sigma}-\gamma^{-2}\zeta^2, 
\end{equation*}
where the time derivative is calculated along the trajectories of \eqref{Nmodes}, \eqref{output}, \eqref{e} with $L$ and $K$ given in \eqref{controlGains} with $X$ and $Y$ replaced by $\tilde X$ and $\tilde Y$. For the zero initial conditions, $\tilde V(0)=0$. Since $\lim_{t\to\infty}\tilde V(t)\ge0$, integrating the above from $0$ to $\infty$, we find 
\begin{equation}\label{HinfSigma}
    \frac{\sigma}{\gamma}\|z^N\|^2+\|u\|^2\le\frac{\|\zeta\|^2}{\gamma^2}+\frac{1}{\gamma\sigma}\|F\|^2, 
\end{equation}
where $\|\cdot\|$ is the norm in $L^2(0,\infty)$. The open-loop system \eqref{Nmodes}, \eqref{output} can be written as 
\begin{equation*}
\begin{aligned}
\dot z^{N}&=Az^{N}+B_1w+Bu,\\
\tilde z&=C_1z^{N}+D_{12}u,\\
y&=Cz^{N}+D_{21}w, 
\end{aligned}
\end{equation*}
where $A$, $B$, and $C$ are from \eqref{notations}, and 
\begin{equation*}
\begin{aligned}
&B_1=\begin{bmatrix}
0_{N\times 1} & \sqrt{\sigma/\gamma}I_{N}
\end{bmatrix},&&D_{21}=\begin{bmatrix}
1&0_{1\times N}
\end{bmatrix},\\
&C_1=\begin{bmatrix}
\sqrt{\sigma/\gamma}I_{N}\\0_{1\times N}
\end{bmatrix},\quad 
&&D_{12}=\begin{bmatrix}
0_{N\times 1}\\1
\end{bmatrix},
\end{aligned}
\end{equation*}
$w=\operatorname{col}\{\zeta,\sqrt{\gamma/\sigma}F\}$, and $\tilde z$ is the controlled output. Relation \eqref{HinfSigma} implies
\begin{equation*}
    \|\tilde z\|\le\gamma^{-1}\|w\|. 
\end{equation*}
That is, there is a linear controller leading to the $L^2$ gain not greater than $\gamma^{-1}$. By \cite[Theorem~8.3.2]{Green2012}, there must exist non-negative stabilizing solutions of the algebraic Riccati equations \eqref{ARE} such that $\rho(XZ)<\gamma^{-2}$. Clearly, these solutions are positive definite since $\frac{\sigma}{\gamma}I_N>0$ in \eqref{ARE}. 
\section{Proof of Proposition~\ref{prop:N}}\label{prop:N:proof}\color{black}
The conditions of Theorem~\ref{th:stability} guarantee that both \eqref{Vzdot} and \eqref{Vedot} are not positive, that is, 
\begin{equation}\label{VzVe}
    \textstyle\dot V_z+\dot V_e\le-u^2-\frac{\sigma}{\gamma}|z^N|^2+\frac1{\gamma\sigma}|F|^2+\frac{\zeta^2}{\gamma^2}. 
\end{equation}
For $V_N=\gamma^{-1}z_N^2$ and $\mu_N>0$, Young's inequality gives 
\begin{equation}\label{VN}
\begin{aligned}
    \dot{V}_N&\textstyle\!\stackrel{\eqref{residue}}{=}\!2\gamma^{-1} z_N\left(-\left(\lambda_N-q\right) z_N+b_Nu+f_N\right)\\
    &\textstyle\le \left[-\frac{2(\lambda_N-q)}{\gamma}+\frac{b_N^2}{\gamma^2}\mu_N+\frac{\sigma}{\gamma}\right]z_N^2+\frac{u^2}{\mu_N}+\frac{f_N^2}{\gamma\sigma}\\
    &\textstyle\!\stackrel{\eqref{mu_n}}{=}\!\left[-\frac{2\mu_N}{\pi\gamma^2}-\frac{\sigma}{\gamma}\right]z_N^2+\frac{u^2}{\mu_N}+\frac{f_N^2}{\gamma\sigma}.
\end{aligned}
\end{equation}
By increasing the number of considered modes from $N$ to $N+1$, the input-to-residue $L^2$ gain changes to (see \eqref{gamma:series}) 
\begin{equation}\label{gammabar}
\textstyle\bar\gamma=\frac{2}{\pi}\sum_{n=N+1}^\infty\frac{1}{\lambda_n-q-\sigma}\stackrel{\eqref{mu_n}}{=}\gamma\sum_{n=N+1}^\infty\mu_n^{-1}. 
\end{equation}

Consider $\bar V=\frac{\gamma}{\bar\gamma}(V_z+V_e+V_N)$. Multiplying \eqref{VzVe} and \eqref{VN} by $\frac{\gamma}{\bar\gamma}$, we find 
\begin{equation*}
\textstyle\dot{\bar V}\le-\frac{\gamma}{\bar\gamma}(1-\mu_N^{-1})u^2-\frac{\sigma}{\bar\gamma}\left|z^{N+1}\right|^2-\frac{2\mu_Nz_N^2}{\pi\gamma\bar\gamma}+\frac{|\bar F|^2}{\bar\gamma\sigma}+\frac{\zeta^2}{\gamma\bar\gamma},  
\end{equation*}
where $\bar F=[f_0,\ldots,f_N]^\top$. Recall that $\gamma$ was selected so that $\sum_{n=N}^\infty\mu_n^{-1}=1$. Thus, the coefficient in front of $u^2$ is 
\begin{equation*}
\textstyle\frac{\gamma}{\bar\gamma}(1-\mu_N^{-1})=\frac{\gamma}{\bar\gamma}\sum_{n=N+1}^\infty\mu_n^{-1}\stackrel{\eqref{gammabar}}{=}1. 
\end{equation*}
Furthermore, \eqref{gammabar} and Young's inequality imply 
\begin{equation*}
\begin{aligned}
\textstyle\frac{\zeta^2}{\gamma\bar\gamma}-\frac{2\mu_Nz_N^2}{\pi\gamma\bar\gamma}&=\bar\gamma^{-2}\left[(1-\mu_N^{-1})\zeta^2+(1-\mu_N)c_N^2z_N^2\right]\\
&\le\bar\gamma^{-2}\left[\zeta^2-2\zeta c_Nz_N+c_N^2z_N^2\right]=\bar\gamma^{-2}\bar\zeta^2, 
\end{aligned}
\end{equation*}
where 
\begin{equation*}
\textstyle\bar\zeta=\sum_{n=N+1}^\infty c_nz_n=\zeta-c_Nz_N. 
\end{equation*}
Combining the last two bounds, we obtain (cf.~\eqref{VzVe})
\begin{equation*}
\textstyle\dot{\bar V}\le-u^2-\frac{\sigma}{\bar\gamma}\left|z^{N+1}\right|^2+\frac{1}{\bar\gamma\sigma}|\bar F|^2+\frac{\bar\zeta^2}{\bar\gamma^2}.
\end{equation*}
For the zero initial conditions, $\bar V(0)=0$. Since $\lim_{t\to\infty}\bar V(t)\ge0$, integrating the above from $0$ to $\infty$, we find 
\begin{equation}\label{HinfN+1}
\textstyle\frac{\sigma}{\bar\gamma}\|z^{N+1}\|^2+\|u\|^2\le\frac{\|\bar\zeta\|^2}{\bar\gamma^2}+\frac{1}{\bar\gamma\sigma}\|\bar F\|^2, 
\end{equation}
where $\|\cdot\|$ is the norm in $L^2(0,\infty)$. The open-loop system \eqref{Nmodes}, \eqref{output} with $N$ replaced by $N+1$ can be written as 
\begin{equation*}
\begin{aligned}
\dot z^{N+1}&=\bar Az^{N+1}+B_1w+B_2u,\\
\tilde z&=C_1z^{N+1}+D_{12}u,\\
y&=C_2z^{N+1}+D_{21}w, 
\end{aligned}
\end{equation*}
where 
\begin{equation*}\textcolor{change}{
\begin{aligned}
&\bar A=\operatorname{diag}\{q-\lambda_0,\ldots,q-\lambda_{N}\},\\
&B_1=\begin{bmatrix}
0_{(N+1)\times 1} & \sqrt{\sigma/\bar\gamma}I_{N+1}
\end{bmatrix},&&B_2=[b_0,b_1,\ldots,b_{N}]^\top,\\
&C_1=\begin{bmatrix}
\sqrt{\sigma/\bar\gamma}I_{N+1}\\0_{1\times(N+1)}
\end{bmatrix},\quad 
&&D_{12}=\begin{bmatrix}
0_{(N+1)\times 1}\\1
\end{bmatrix},\\
&C_2=[c_0,c_1,\ldots,c_{N}],\quad &&D_{21}=\begin{bmatrix}
1&0_{1\times(N+1)}
\end{bmatrix},
\end{aligned}}
\end{equation*}
$w=\operatorname{col}\{\bar\zeta,\sqrt{\bar\gamma/\sigma}\bar F\}$, and $\tilde z$ is the controlled output. Relation \eqref{HinfN+1} takes the form 
\begin{equation*}
    \|\tilde z\|\le\bar\gamma^{-1}\|w\|. 
\end{equation*}
That is, there is a linear controller leading to the $L^2$ gain not greater than $\bar\gamma^{-1}$. By \cite[Theorem~8.3.2]{Green2012}, there must exist non-negative stabilizing solutions of the algebraic Riccati equations \eqref{ARE} such that $\rho(XZ)<\bar\gamma^{-2}$. Clearly, these solutions are positive definite since $\frac{\sigma}{\bar\gamma}I_{N+1}>0$ in \eqref{ARE}. 
\section{Proof of Proposition~\ref{prop:h}}\label{prop:h:proof}
Consider $\tilde\Psi$ comprising the first four block rows and columns of $\Psi$ defined in \eqref{Psi} $q$ replaced by $\tilde q$ in $A$. If $\sigma=0$, eliminate the third block column and row from $\Psi$ and $\tilde\Psi$. The Schur complement lemma implies that $\tilde\Psi<0$ is equivalent to 
\begin{equation}\label{PsiSchur}
\left[\begin{smallmatrix}
\Psi_{11}+\gamma\sigma P_z^2 & \Psi_{12}-\gamma\sigma P_zP_e \\
* & \Psi_{22}+\gamma^2P_eLL^\top P_e+\gamma\sigma P_e^2
\end{smallmatrix}\right]<0. 
\end{equation}
 Let $X>0$ and $Z>0$ be the stabilizing solutions of \eqref{ARE} and $Y=\gamma^{-2}Z^{-1}-X$. Recall that $Y>0$ by \eqref{spectralCondition}. Taking $P_z=X$ and $P_e=Y$, and substituting the expressions for $K$ and $L$ from \eqref{controlGains}, we find that \eqref{PsiSchur} is equivalent to 
\begin{equation*}
\begin{aligned}
&X\tilde A+\tilde A^{\top}X-X\left(BB^{\top}-\gamma\sigma I_N\right)X+\frac{\sigma}{\gamma}I_N<0,\\
&Y(\tilde A+\gamma\sigma X)+(\tilde A+\gamma\sigma X)^{\top} Y+XBB^{\top}X\notag \\
&\hspace{4.9cm}{}+\gamma\sigma Y^2-\gamma^{-2}C^\top C<0, 
\end{aligned}
\end{equation*}
where $\tilde A$ is obtained from $A$ by replacing $q$ with $\tilde q$. Note that the off-diagonal block in \eqref{PsiSchur} becomes zero. 
Substituting $\tilde q=q-\varepsilon>0$ with $\varepsilon>0$, we obtain 
\begin{equation*}
\begin{aligned}
&X\tilde A+\tilde A^{\top}X-X\left(BB^{\top}-\gamma\sigma I_N\right)X+\frac{\sigma}{\gamma}I_N\\
&=\eqref{ARE4X}-2\varepsilon X=-2\varepsilon X<0, \\
&Y(\tilde A+\gamma\sigma X)+(\tilde A+\gamma\sigma X)^{\top} Y+XBB^{\top}X\!\\
&\hspace{5.2cm}+\gamma\sigma Y^2-\gamma^{-2}C^\top C\\
&=\gamma^{-2}Z^{-1}\times\eqref{ARE4Z}\times Z^{-1}-\eqref{ARE4X}-2\varepsilon Y=-2\varepsilon Y<0. 
\end{aligned}
\end{equation*}
This implies \eqref{PsiSchur} and, therefore, $\tilde\Psi<0$. Then, $\Psi<0$ holds for $h=0$ and $W_z=W_e=\alpha I_N$ with a large enough $\alpha>0$. By continuity, $\Psi<0$ for a small enough $h>0$.
\end{document}